# Average velocity and effective diffusion of a Brownian particle driven by a constant force over a static periodic potential

## Hongyun Wang

Department of Applied Mathematics and Statistics, Jack Baskin School of Engineering,
University of California, Santa Cruz, CA 95064, USA

Email: hongwang@soe.ucsc.edu

## Abstract

In this manuscript, we consider the case where a Brownian particle is subject to a static periodic potential and is driven by a constant force. We derive analytic formulas for the average velocity and the effective diffusion.

**Mathematical formulation**

We consider the one-dimensional stochastic motion of a particle subject to a periodic potential and a constant driving force. In particular, we are interested in the effective diffusion coefficient of the particle, and how the periodic potential and the driving force affect the effective diffusion coefficient.

The probability density of the particle is governed by the Fokker-Planck equation:

$$\frac{\partial \rho}{\partial t} = D \frac{\partial}{\partial x} \left( \frac{(\phi' - f)}{k_B T} \rho + \frac{\partial \rho}{\partial x} \right) \tag{1}$$

where $\rho(x,t)$ is the probability density that the particle is at position $x$ at time $t$, $D$ is the diffusion constant of the particle, $\phi(x)$ is the periodic potential with period $= L$, $f$ is the constant driving force, $k_B$ is the Boltzmann constant, and $T$ is the absolute temperature.

The probability density $\rho(x,t)$ satisfies the normalizing condition

$$\int_{-\infty}^{\infty} \rho(x,t)dx = 1$$

Consider the transform $x = L\tilde{x}$, $t = \dfrac{L^2}{D}\tilde{t}$

The function $\tilde{\rho}(\tilde{x},\tilde{t}) = L\rho\left(L\tilde{x}, \dfrac{L^2\tilde{t}}{D}\right)$ satisfies

$$\frac{\partial \tilde{\rho}}{\partial \tilde{t}} = \frac{\partial}{\partial \tilde{x}}\left((\tilde{\phi}' - \tilde{f})\tilde{\rho} + \frac{\partial \tilde{\rho}}{\partial \tilde{x}}\right), \qquad \int_{-\infty}^{\infty} \tilde{\rho}(\tilde{x},\tilde{t})d\tilde{x} = 1 \qquad (2)$$

where $\tilde{\phi}(\tilde{x}) = \dfrac{\phi(L\tilde{x})}{k_B T}$ is a periodic function of $\tilde{x}$ with period $= 1$ and $\tilde{f} = \dfrac{fL}{k_B T}$ is a constant.

Thus, we only need to consider equation (1) for

$$L = 1, \quad D = 1, \quad k_B T = 1$$

**The effective diffusion coefficient**

The mean and the variance of the particle position at time $t$ are

$$\langle X(t) \rangle = \int_{-\infty}^{\infty} x\,\rho(x,t)\,dx$$

$$\text{var}\{X(t)\} = \left\langle \left(X(t) - \langle X(t)\rangle\right)^2 \right\rangle = \int_{-\infty}^{\infty} \left(x - \langle X(t)\rangle\right)^2 \rho(x,t)\,dx$$

It can be shown that for large $t$, both the mean and the variance increase approximately linearly. The average velocity and the effective diffusion coefficient of the particle are defined as

$$V \stackrel{def}{=} \lim_{t \to \infty} \frac{\langle X(t)\rangle}{t}$$

$$D_{\text{eff}} \stackrel{def}{=} \lim_{t \to \infty} \frac{\text{var}\{X(t)\}}{2t}$$

The effective drag coefficient can be defined as

$$\zeta_{\text{eff}} \stackrel{def}{=} \frac{f}{V}$$

Let us introduce functions

$$\rho_k(x,t) = \sum_{j=-\infty}^{\infty} j^k \rho(j+x,t) \tag{3}$$

In appendix A, we derive that

$$\lim_{t\to\infty} \frac{1}{t}\int_0^1 \rho_1(x,t)\,dx = V \tag{4}$$

$$\lim_{t\to\infty} \frac{1}{2t}\left[\int_0^1 \rho_2(x,t)\,dx - \left(\int_0^1 \rho_1(x,t)\,dx\right)^2\right] = D_{\text{eff}} \tag{5}$$

Below we give the formulas for the average velocity and effective diffusion. The detailed derivations are given in Appendices B, C and D.

The steady state of $\rho_0(x,t)$ is

$$u_0(x) = \frac{w_0(x)}{\int_0^1 w_0(x)\,dx}$$

where

$$w_0(x) = \int_0^1 \exp(\phi(x+s) - \phi(x) - f s)\,ds$$

The average velocity $V$, the effective diffusion $D_{\text{eff}}$, and the effective drag coefficient $\zeta_{\text{eff}}$ are given by (see Appendices B, C and D for derivations)

$$V = \frac{1 - \exp(-f)}{\int_0^1 w_0(x)\,dx},$$

$$D_{\text{eff}} = \frac{\int_0^1 (w_0(x))^2 \int_0^1 \exp(\phi(x) - \phi(x-s) - f s)\,ds\,dx}{\left(\int_0^1 w_0(x)\,dx\right)^3}$$

$$\zeta_{\text{eff}} = \frac{f}{1 - \exp(-f)} \int_0^1 w_0(x)\,dx$$

**Results and discussion**

Below we discuss the behaviors of effective diffusion and effective drag coefficient for several cases (the detailed derivation is in Appendix E).

<u>Case A:</u>   $f$ is small

$$V = \frac{f}{a_0}\left(1 + f\left(\frac{a_1}{a_0} - \frac{1}{2}\right) + O(f^2)\right)$$

$$\zeta_{\text{eff}} = a_0\left(1 + f\left(\frac{1}{2} - \frac{a_1}{a_0}\right) + O(f^2)\right)$$

$$D_{\text{eff}} = \frac{1}{a_0}\left(1 + f\left(\frac{a_1}{a_0} - \frac{1}{2}\right) + O(f^2)\right)$$

$$\zeta_{\text{eff}} D_{\text{eff}} = \left(1 + O(f^2)\right)$$

where

$$a_0 = \left(\int_0^1 \exp(-\phi(x))dx\right)\left(\int_0^1 \exp(\phi(s))ds\right)$$

$$a_1 = \int_0^1\int_0^1 \exp(-\phi(x) + \phi(x+s))s\,ds\,dx$$

If $\dfrac{a_1}{a_0} \neq \dfrac{1}{2}$ (this happens when the potential is asymmetric), then the minimum of $D_{\text{eff}}$ is attained at non-zero value of $f$, and

$$\min D_{\text{eff}} < \frac{1}{a_0}$$

<u>Case B:</u>   $f$ is large

$$V = f\left[1 - \frac{1}{f^2}\int_0^1 (\phi'(x))^2 dx + O\left(\frac{1}{f^3}\right)\right]$$

$$\zeta_{\text{eff}} = 1 + \frac{1}{f^2}\int_0^1 (\phi'(x))^2 dx + O\left(\frac{1}{f^3}\right)$$

$$D_{eff} = 1 + \frac{3}{f^2}\int_0^1 (\phi'(x))^2 dx + O\left(\frac{1}{f^3}\right)$$

$$\zeta_{eff} D_{eff} = 1 + \frac{4}{f^2}\int_0^1 (\phi'(x))^2 dx + O\left(\frac{1}{f^3}\right)$$

It is clear that for large $f$

$$\zeta_{eff} > \zeta, \quad \text{and} \quad D_{eff} > D$$

In both Case A and Case B, we see that the non-dimensional Einstein relation $\zeta_{eff} D_{eff} = 1$ is approximately satisfied.

<u>Case A:</u> $\quad \zeta_{eff} D_{eff} = 1 + O(f^2) \quad$ for small $f$

<u>Case B:</u> $\quad \zeta_{eff} D_{eff} = 1 + O\left(\frac{1}{f^2}\right) \quad$ for large $f$

<u>Case C:</u> $\quad \phi(x) = \begin{cases} -A & 0 \leq x < 0.5 \\ A & 0.5 \leq x < 1 \end{cases}$

## Appendix A

In this appendix, we derive equations (4) and (5):

$$\lim_{t\to\infty} \frac{1}{t}\int_0^1 \rho_1(x,t)\,dx = V$$

$$\lim_{t\to\infty} \frac{1}{2t}\left[\int_0^1 \rho_2(x,t)\,dx - \left(\int_0^1 \rho_1(x,t)\,dx\right)^2\right] = D_{eff}$$

We first write $\int_0^1 \rho_1(x,t)\,dx$ as

$$\int_0^1 \rho_1(x,t)\,dx = \int_0^1 \sum_{j=-\infty}^{\infty}(j+x)\rho(j+x,t)\,dx - \int_0^1 x\rho_0(x,t)\,dx \equiv \langle X(t)\rangle - T_1$$

$T_1 = \int_0^1 x\rho_0(x,t)\,dx$ satisfies $0 \leq T_1 \leq 1$.

$$\Longrightarrow \quad \lim_{t\to\infty} \frac{1}{t}\int_0^1 \rho_1(x,t)\,dx = \lim_{t\to\infty} \frac{1}{t}\big[\langle X(t)\rangle - T_1\big] = V$$

Then we have

$$\int_0^1 \rho_2(x,t)\,dx - \left(\int_0^1 \rho_1(x,t)\,dx\right)^2 = \int_0^1 \sum_{j=-\infty}^{\infty}\left[j - \int_0^1 \rho_1(s,t)\,ds\right]^2 \rho(j+x,t)\,dx$$

$$= \int_0^1 \sum_{j=-\infty}^{\infty} \big\{[j+x-\langle X(t)\rangle] + [T_1 - x]\big\}^2 \rho(j+x,t)\,dx$$

$$= \text{var}\{X(t)\} + 2\int_0^1 \sum_{j=-\infty}^{\infty} [j+x-\langle X(t)\rangle]\cdot[T_1-x]\rho(j+x,t)\,dx$$

$$+ \int_0^1 \sum_{j=-\infty}^{\infty} [T_1-x]^2 \rho(j+x,t)\,dx \equiv \text{var}\{X(t)\} + T_2 + T_3$$

Using the Cauchy-Schwarz inequality, we see that $T_2$ is bounded by

$$|T_2| \le 2\int_0^1 \sum_{j=-\infty}^{\infty} |j+x-\langle X(t)\rangle|\rho(x+j,t)\,dx$$

$$\le 2\sqrt{\int_0^1 \sum_{j=-\infty}^{\infty} [j+x-\langle X(t)\rangle]^2 \rho(x+j,t)\,dx} = 2\sqrt{\text{var}\{X(t)\}}$$

which leads to

$$\lim_{t\to\infty} \frac{1}{2t}\left[\int_0^1 \rho_2(x,t)\,dx - \left(\int_0^1 \rho_1(x,t)\,dx\right)^2\right] = \lim_{t\to\infty} \frac{1}{2t}\big[\text{var}\{X(t)\} + T_2 + T_3\big] = D_{\text{eff}}$$

### Appendix B

In this appendix, we derive the results for $\rho_0(x,t)$:

- $\rho_0(x,t)$ converges to a steady state $u_0(x)$ given by

$$u_0(x) = \frac{1}{M_0} w_0(x)$$

where

$$w_0(x) = \int_0^1 \exp(\phi(x+s) - \phi(x) - fs)ds, \qquad M_0 = \int_0^1 w_0(x)dx$$

- The steady state flux of $u_0(x)$ is

$$J_0 = -((\phi' - f)u_0(x) + u_0'(x)) = \frac{1}{M_0}(1 - \exp(-f))$$

The boundary condition for $\rho_0(x,t)$ is $\rho_0(x+1,t) = \rho_0(x,t)$.

So $\rho_0(x,t)$ is governed by

$$\begin{cases} \dfrac{\partial \rho_0}{\partial t} = \dfrac{\partial}{\partial x}\left((\phi' - f)\rho_0 + \dfrac{\partial \rho_0}{\partial x}\right) \\ \rho_0(x+1,t) = \rho_0(x,t), \qquad \int_0^1 \rho_0(x,t)dx = 1 \end{cases} \qquad (6)$$

Let $u_0(x)$ be the steady state solution of (6). $u_0(x)$ satisfies

$$(\phi' - f)u_0 + u_0' = -J_0$$

$$\Longrightarrow \quad [\exp(\phi(x) - fx)u_0(x)]' = -J_0 \exp(\phi(x) - fx)$$

$$\Longrightarrow \quad u_0(x) = \exp(-\phi(x) + fx)\left[c_1 + \frac{J_0}{1 - \exp(-f)}\int_x^{x+1} \exp(\phi(s) - fs)ds\right]$$

The boundary condition: $u_0(x+1) = u_0(x)$ implies $c_1 = 0$. Thus, we get

$$u_0(x) = \frac{J_0}{1 - \exp(-f)} w_0(x)$$

where

$$w_0(x) = \int_0^1 \exp(\phi(x+s) - \phi(x) - fs)ds.$$

Applying the normalizing condition: $\int_0^1 u_0(x)dx = 1$, we have

$$J_0 = \frac{1}{M_0}(1 - \exp(-f)),$$

where

$$M_0 = \int_0^1 w_0(x)dx$$

Now we show that $\lim_{t\to\infty} \rho_0(x,t) = u_0(x)$

We consider $r_0(x,t) = \frac{\rho_0(x,t)}{u_0(x)} - 1$. $r_0(x,t)$ satisfies

$$\begin{cases} u_0 \frac{\partial r_0}{\partial t} = \frac{\partial}{\partial x}\left(-J_0 r_0 + u_0 \frac{\partial r_0}{\partial x}\right) \\ r_0(x+1,t) = r_0(x,t), \quad \int_0^1 u_0(x) r_0(x,t) dx = 0 \end{cases} \quad (7)$$

To show $\lim_{t\to\infty} r_0(x,t) = 0$, we consider the time evolution of $\int_0^1 u_0(x)(r_0(x,t))^2 dx$.

$$\frac{d}{dt}\int_0^1 u_0 r_0^2 dx = 2\int_0^1 u_0 r_0 \frac{\partial r_0}{\partial t} dx = 2\int_0^1 r_0 \frac{\partial}{\partial x}\left(-J_0 r_0 + u_0 \frac{\partial r_0}{\partial x}\right) dx$$

$$= -2\int_0^1 \frac{\partial r_0}{\partial x}\left(-J_0 r_0 + u_0 \frac{\partial r_0}{\partial x}\right) dx = -2\int_0^1 u_0\left(\frac{\partial r_0(x,t)}{\partial x}\right)^2 dx \leq 0$$

The quantity $\int_0^1 u_0 r_0^2 dx$ is non-negative and non-increasing. Hence, we obtain

$$\lim_{t\to\infty} \int_0^1 u_0\left(\frac{\partial r_0(x,t)}{\partial x}\right)^2 dx = 0$$

$$\Longrightarrow \quad \lim_{t\to\infty} \int_0^1 \left|\frac{\partial r_0(x,t)}{\partial x}\right| dx = 0$$

Using the constraint $\int_0^1 u_0(x) r_0(x,t) dx = 0$, we write $r_0(x,t)$ as

$$r_0(x,t) = -\int_0^1 \left(\int_0^\xi \frac{\partial r_0(s,t)}{\partial x} ds\right) u_0(\xi) d\xi + \int_0^x \frac{\partial r_0(s,t)}{\partial x} ds$$

Thus, we have

$$\lim_{t\to\infty} r_0(x,t) = 0$$

==> $$\lim_{t\to\infty} \rho_0(x,t) = u_0(x)$$

For simplicity, we start with the initial probability density

$$\rho(x,0) = \begin{cases} u_0(x) & 0 \le x \le 1 \\ 0 & \text{otherwise} \end{cases}$$

With this initial probability density, we have

$$\rho_0(x,t) = u_0(x) \quad \text{and} \quad \rho_1(x,0) \equiv 0$$

**Appendix C**

In this appendix, we derive the results for $\rho_1(x,t)$:

- $p_1(x,t) = \rho_1(x,t) - u_0(x)J_0 t$ converges to a steady state $u_1(x)$ given by

$$u_1(x) = \left(\frac{M_1}{M_0^3 J_0} + \frac{1}{2}\right) u_0(x) - u_0(x) \int_0^x u_0(s)ds - \frac{w_1(x)}{M_0^3 J_0}$$

where

$$w_1(x) = \int_0^1 (w_0(x+s))^2 \exp(\phi(x+s) - \phi(x) - fs)ds,$$

$$M_1 = \int_0^1 w_1(x)dx$$

- The steady state flux of $u_1(x)$ is

$$-((\phi' - f)u_1 + u_1') = \frac{M_1}{M_0^3} + \frac{J_0}{2} - J_0 \int_0^x u_0(s)ds$$

The boundary condition for $\rho_1(x,t)$ is $\rho_1(x+1,t) = \rho_1(x,t) - u_0(x)$.

So $\rho_1(x,t)$ is governed by

$$\begin{cases} \dfrac{\partial \rho_1}{\partial t} = \dfrac{\partial}{\partial x}\left((\phi' - f)\rho_1 + \dfrac{\partial \rho_1}{\partial x}\right) \\ \rho_1(x+1,t) = \rho_1(x,t) - u_0(x) \end{cases} \qquad (8)$$

Integrating (8), we obtain

$$\frac{d}{dt}\int_0^1 \rho_1 dx = \left((\phi' - f)\rho_1 + \frac{\partial \rho_1}{\partial x}\right)\bigg|_{x=0}^{x=1} = -((\phi' - g)u_0(0) + u_0'(0)) = J_0 \qquad (9)$$

This shows that $\rho_1(x,t)$ does not converge to a steady state.

$p_1(x,t) = \rho_1(x,t) - u_0(x)J_0 t$ is governed by

$$\begin{cases} \frac{\partial p_1}{\partial t} + J_0 u_0 = \frac{\partial}{\partial x}\left((\phi' - f)p_1 + \frac{\partial p_1}{\partial x}\right) \\ p_1(x+1,t) = p_1(x,t) - u_0(x), \qquad \int_0^1 p_1(x,t)dx = 0 \end{cases} \qquad (10)$$

Let $u_1(x)$ be the stead state solution of (10). $u_1(x)$ satisfies

$$\begin{cases} J_0 u_0 = \left[(\phi' - f)u_1 + u_1'\right]' \\ u_1(x+1) = u_1(x) - u_0(x), \qquad \int_0^1 u_1(x)dx = 0 \end{cases}$$

$$\Longrightarrow \quad (\phi' - f)u_1 + u_1' = -c_1 + J_0 \int_0^x u_0(s)ds$$

$$\Longrightarrow \quad \left[\exp(\phi(x) - fx)u_1(x)\right]' = \exp(\phi(x) - fx)\left(-c_1 + J_0 \int_0^x u_0(s)ds\right)$$

Using $\exp(\phi(x) - fx) = \frac{-1}{J_0}\left[\exp(\phi(x) - fx)u_0(x)\right]'$, we obtain

$$u_1(x) = \frac{c_1}{J_0}u_0(x) - \exp(-\phi(x) + fx)\int \left[\exp(\phi(x) - fx)u_0(x)\right]' \int_0^x u_0(s)ds\, dx$$

$$= \frac{c_1}{J_0}u_0(x) - u_0(x)\int_0^x u_0(s)ds + \exp(-\phi(x) + fx)\int \exp(\phi(x) - fx)(u_0(x))^2 dx$$

$$= \frac{c_1}{J_0}u_0(x) - u_0(x)\int_0^x u_0(s)ds - \frac{\exp(-\phi(x) + fx)}{1 - \exp(-f)}\left(\int_x^{x+1} \exp(\phi(s) - fs)(u_0(s))^2 ds + c_2\right)$$

The condition $u_1(x+1) = u_1(x) - u_0(x)$ implies $c_2 = 0$. Thus, we obtain

$$u_1(x) = \frac{c_1}{J_0} u_0(x) - u_0(x) \int_0^x u_0(s) ds - \frac{1}{M_0^3 J_0} w_1(x)$$

where

$$w_1(x) = \int_0^1 (w_0(x+s))^2 \exp(\phi(x+s) - \phi(x) - fs) ds$$

The condition $\int_0^1 u_1(x) dx = 0$ yields $c_1 = \frac{M_1}{M_0^3} + \frac{J_0}{2}$ where $M_1 = \int_0^1 w_1(x) dx$

Thus, the steady state flux of $u_1(x)$ is

$$-((\phi' - f) u_1 + u_1') = \frac{M_1}{M_0^3} + \frac{J_0}{2} - J_0 \int_0^x u_0(s) ds$$

Now we show that $\lim_{t \to \infty} p_1(x,t) = u_1(x)$

We consider $r_1(x,t) = \frac{p_1(x,t)}{u_0(x)} - \frac{u_1(x)}{u_0(x)}$. $r_1(x,t)$ is governed by

$$\begin{cases} u_0 \dfrac{\partial r_1}{\partial t} = \dfrac{\partial}{\partial x}\left(-J_0 r_1 + u_0 \dfrac{\partial r_1}{\partial x}\right) \\ r_1(x+1, t) = r_1(x,t), \quad \int_0^1 u_0(x) r_1(x,t) dx = 0 \end{cases}$$

This is the same system as (7)

Thus, we have

$$\lim_{t \to \infty} r_1(x,t) = 0 \quad \Longrightarrow \quad \lim_{t \to \infty} p_1(x,t) = u_1(x)$$

**Appendix D**

In this appendix, we derive the effective diffusion.

The boundary condition for $p_2(x,t)$ is $p_2(x+1,t) = p_2(x,t) - 2p_1(x,t) + u_0(x)$.

So $p_2(x,t)$ is governed by

$$\begin{cases} \dfrac{\partial \rho_2}{\partial t} = \dfrac{\partial}{\partial x}\left((\phi' - f)\rho_2 + \dfrac{\partial \rho_2}{\partial x}\right) \\ \rho_2(x+1,t) = \rho_2(x,t) - 2\rho_1(x,t) + u_0(x) \end{cases} \quad (11)$$

Integrating (11) and using $\rho_2(x+1,t) = \rho_2(x,t) - 2\rho_1(x,t) + u_0(x)$ and $\rho_1(x,t) = p_1(x,t) + u_0(x)J_0 t$, we obtain

$$\dfrac{d}{dt}\int_0^1 \rho_2 dx = \left((\phi'-f)\rho_2 + \dfrac{\partial \rho_2}{\partial x}\right)\Big|_{x=0}^{x=1} = -2\left((\phi'-f)\rho_1 + \dfrac{\partial \rho_1}{\partial x}\right)\Big|_{x=0} - J_0$$

$$= -2\left((\phi'-f)p_1 + \dfrac{\partial p_1}{\partial x}\right)\Big|_{x=0} + (2J_0^2 t - J_0)$$

Using the steady state flux of $u_1(x)$, we get

$$\lim_{t\to\infty} \dfrac{d}{dt}\int_0^1 \rho_2 dx = -2\left((\phi'-f)u_1 + \dfrac{\partial u_1}{\partial x}\right)\Big|_{x=0} + (2J_0^2 t - J_0) = 2\dfrac{M_1}{M_0^3} + 2J_0^2 t$$

$$\implies \lim_{t\to\infty} \dfrac{d}{dt}\left[\int_0^1 \rho_2 dx - \left(\int_0^1 \rho_1 dx\right)^2\right] = \lim_{t\to\infty} \dfrac{d}{dt}\int_0^1 \rho_2 dx - 2J_0^2 t = 2\dfrac{M_1}{M_0^3}$$

$$\implies D_{\text{eff}} = \lim_{t\to\infty} \dfrac{1}{2t}\left[\int_0^1 \rho_2 dx - \left(\int_0^1 \rho_1 dx\right)^2\right] = \dfrac{M_1}{M_0^3}$$

The average velocity is

$$V = J_0 = \dfrac{1 - \exp(-f)}{M_0}$$

The effective drag coefficient is

$$\zeta_{\text{eff}} = \dfrac{f}{V} = \dfrac{f M_0}{1 - \exp(-f)}$$

Notice that

$$M_1 = \int_0^1 w_1(x)dx = \int_0^1\int_0^1 (w_0(x+s))^2 \exp(\phi(x+s) - \phi(x) - fs)dsdx$$

$$= \int_0^1 \int_0^1 (w_0(x))^2 \exp(\phi(x) - \phi(x-s) - fs)\,ds\,dx$$

**Appendix E**

In this appendix, we derive the asymptotic behaviors of the effective diffusion and the effective drag coefficient for several cases

Case A:  $f$ is small

$$1 - \exp(-f) = f\left[1 - \frac{f}{2} + O(f^2)\right]$$

$$\frac{f}{1 - \exp(-f)} = 1 + \frac{f}{2} + O(f^2)$$

$$w_0(x) = \int_0^1 \exp(\phi(x+s) - \phi(x) - fs)\,ds$$

$$= \exp(-\phi(x)) \int_0^1 \exp(\phi(x+s))[1 - fs + O(f^2)]\,ds$$

$$\implies \int_0^1 w_0(x)\,dx = a_0\left[1 - f\frac{a_1}{a_0} + O(f^2)\right]$$

where $a_0 = \left(\int_0^1 \exp(-\phi(x))\,dx\right)\left(\int_0^1 \exp(\phi(s))\,ds\right)$

$$a_1 = \int_0^1 \int_0^1 \exp(-\phi(x) + \phi(x+s))s\,ds\,dx$$

$$\implies V = \frac{1 - \exp(-f)}{\int_0^1 w_0(x)\,dx} = \frac{f}{a_0}\left[1 + f\left(\frac{a_1}{a_0} - \frac{1}{2}\right) + O(f^2)\right]$$

$$\zeta_{\text{eff}} = \frac{f \int_0^1 w_0(x)\,dx}{1 - \exp(-f)} = a_0\left[1 + f\left(\frac{1}{2} - \frac{a_1}{a_0}\right) + O(f^2)\right]$$

$$\int_0^1 (w_0(x))^2 \int_0^1 \exp(\phi(x) - \phi(x-s) - fs) ds\, dx$$

$$= \int_0^1 e^{-\phi(x)} \left\{ \int_0^1 e^{\phi(s)} ds - f \int_0^1 e^{-\phi(x-s)} s\, ds \right\} \left[ \left( \int_0^1 e^{\phi(s)} ds \right)^2 - 2f \left( \int_0^1 e^{\phi(s)} ds \int_0^1 e^{\phi(x+s)} s\, ds \right) \right] dx$$

$$= a_0^2 - f\left(2a_0 a_1 + \frac{1}{2} a_0^2\right) + O(f^2) = a_0^2 \left[ 1 - f\left(\frac{2a_1}{a_0} + \frac{1}{2}\right) + O(f^2) \right]$$

$$\left( \int_0^1 w_0(x) dx \right)^{-3} = \frac{1}{a_0^3} \left[ 1 + 3f \frac{a_1}{a_0} + O(f^2) \right]$$

==> $\quad D_{\text{eff}} = \dfrac{\int_0^1 (w_0(x))^2 \int_0^1 \exp(\phi(x) - \phi(x-s) - fs) ds\, dx}{\left( \int_0^1 w_0(x) dx \right)^3} = \dfrac{1}{a_0} \left[ 1 + f\left(\dfrac{a_1}{a_0} - \dfrac{1}{2}\right) + O(f^2) \right]$

$$\zeta_{\text{eff}} D_{\text{eff}} = 1 + O(f^2)$$

If the potential $\phi(x)$ is an even function, then

$$a_1 = \int_0^1 \int_0^1 \exp(\phi(x+s) - \phi(x)) s\, ds\, dx = \int_0^1 \int_0^1 \exp(\phi(x+1-s) - \phi(x))(1-s) ds\, dx$$

$$= \int_0^1 \int_0^1 \exp(\phi(-x-s) - \phi(-x))(1-s) ds\, dx = \int_0^1 \int_0^1 \exp(\phi(x+s) - \phi(x))(1-s) ds\, dx$$

$$= a_0 - a_1$$

==> $\quad \dfrac{a_1}{a_0} = \dfrac{1}{2}$

In this case, we have

$$\zeta_{\text{eff}} = a_0 \left[ 1 + O(f^2) \right]$$

$$D_{\text{eff}} = \frac{1}{a_0} \left[ 1 + O(f^2) \right]$$

<u>Case B:</u>  $f$ is large

$$f w_0(x) = f \int_0^1 \exp(\phi(x+s) - \phi(x) - fs)ds = \int_0^f \exp\left(\phi\left(x + \frac{s}{f}\right) - \phi(x) - s\right)ds$$

$$= \int_0^f \left[\phi'(x)\frac{s}{f} + \left((\phi'(x))^2 + \phi''(x)\right)\frac{s^2}{2f^2}\right]\exp(-s)ds$$

$$= 1 + \frac{1}{f}\phi'(x) + \frac{1}{f^2}\left((\phi'(x))^2 + \phi''(x)\right) + O\left(\frac{1}{f^3}\right)$$

$$==> \quad f\int_0^1 w_0(x)\,dx = 1 + \frac{1}{f^2}\int_0^1 (\phi'(x))^2\,dx + O\left(\frac{1}{f^3}\right)$$

$$==> \quad V = \frac{1-\exp(-f)}{\int_0^1 w_0(x)\,dx} = f\left[1 - \frac{1}{f^2}\int_0^1 (\phi'(x))^2\,dx + O\left(\frac{1}{f^3}\right)\right]$$

$$\zeta_{eff} = \frac{f\int_0^1 w_0(x)\,dx}{1-\exp(-f)} = 1 + \frac{1}{f^2}\int_0^1 (\phi'(x))^2\,dx + O\left(\frac{1}{f^3}\right)$$

$$f\int_0^1 \exp(\phi(x) - \phi(x-s) - fs)ds = \int_0^f \exp\left(\phi(x) - \phi\left(x - \frac{s}{f}\right) - s\right)ds$$

$$= \int_0^f \left[\phi'(x)\frac{s}{f} + \left((\phi'(x))^2 - \phi''(x)\right)\frac{s^2}{2f^2}\right]\exp(-s)ds$$

$$= 1 + \frac{1}{f}\phi'(x) + \frac{1}{f^2}\left((\phi'(x))^2 - \phi''(x)\right) + O\left(\frac{1}{f^3}\right)$$

$$==> \quad f^3\int_0^1 (w_0(x))^2 \int_0^1 \exp(\phi(x) - \phi(x-s) - fs)\,ds\,dx$$

$$= \int_0^1 \left[1 + \frac{3}{f}\phi'(x) + \frac{1}{f^2}\left(6(\phi'(x))^2 + \phi''(x)\right)\right]dx$$

$$= 1 + \frac{6}{f^2}\int_0^1 (\phi'(x))^2\,dx + O\left(\frac{1}{f^3}\right)$$

$$\implies D_{\text{eff}} = \frac{\int_0^1 (w_0(x))^2 \int_0^1 \exp(\phi(x) - \phi(x-s) - fs)\,ds\,dx}{\left(\int_0^1 w_0(x)\,dx\right)^3}$$

$$= 1 + \frac{3}{f^2}\int_0^1 (\phi'(x))^2\,dx + O\left(\frac{1}{f^3}\right)$$

$$\zeta_{\text{eff}} D_{\text{eff}} = 1 + \frac{4}{f^2}\int_0^1 (\phi'(x))^2\,dx + O\left(\frac{1}{f^3}\right)$$

**Acknowledgements**

This work was partially supported by NSF.